\documentclass[fleqn]{mat01}
\usepackage{times,mathtimy,amssymb}
\begin{document}

\setcounter{page}{103}
\firstpage{103}

\font\zz=msam10 at 10pt
\def\Box{\mbox{\zz{\char'244}}}

\font\zz=msam10 at 10pt
\def\BBox{\mbox{\zz{\char'245}}}

\font\www=mtgub at 10.4pt
\def\piy{\mbox{\www{\char'160}}}
\def\Delt{\mbox{\www{\char'104}}}

\font\bi=tibi at 10.4pt

\newtheorem{theore}{Theorem}
\renewcommand\thetheore{\arabic{section}.\arabic{theore}}
\newtheorem{theor}[theore]{\bf Theorem}
\newtheorem{propo}{\rm PROPOSITION}
\newtheorem{definit}{\rm DEFINITION}
\newtheorem{coro}{\rm COROLLARY}
\newtheorem{rem}[theore]{Remark}
\newtheorem{exampl}[theore]{Example}
\newtheorem{case}{\it Case}

\newtheorem{proo}{Proof}
\renewcommand\theproo{(\arabic{proo})}
\newtheorem{prof}[proo]{Proof of}
\def\exam{\trivlist\item[\hskip\labelsep{\it Example.}]}
\def\pro{\trivlist\item[\hskip\labelsep{PROPOSITION.}]}
\def\theo{\trivlist\item[\hskip\labelsep{\bf Theorem.}]}
\def\lem{\trivlist\item[\hskip\labelsep{\it Lemma.}]}
\def\prof{\trivlist\item[\hskip\labelsep{{\it Proof of} $(*)_{1}$.}]}
\def\prff{\trivlist\item[\hskip\labelsep{{\it Proof of} $(*)_{r} \Rightarrow (*)_{r + 1}$.}]}

\newtheorem{ppro}{Proof}
\renewcommand\theppro{(\arabic{section}) $\Rightarrow$ (\arabic{ppro})}
\newtheorem{proff}[ppro]{Proof of}

\newcommand{\lr}{\longrightarrow}
\newcommand{\ra}{\rightarrow}

\renewcommand{\theequation}{\thesection\arabic{equation}}

\title{A criterion for regular sequences}

\markboth{D~P~Patil, U~Storch and J~St\"{u}ckrad}{A criterion for regular sequences}

\author{D~P~PATIL$^{1}$, U~STORCH$^{2}$ and J~ST\"{U}CKRAD$^{3}$}

\address{$^{1}$Department of Mathematics, Indian Institute of
Science, Bangalore~560~012, India\\
\noindent $^{2}$Fakult\"{a}t f\"{u}r Mathematik, Ruhr Universit\"{a}t Bochum,
D-44780 Bochum, Germany\\
\noindent $^{3}$Fakult\"{a}t f\"{u}r Mathematik und Informatik,
Universit\"{a}t, Leipzig, D-04109 Leipzig, Germany\\
\noindent E-mail: $^{1}$patil@math.iisc.ernet.in; $^{2}$uwe.storch@ruhr-uni-bochum.de;
$^{3}$stueckrad@mathematik.uni-leipzig.de}

\volume{114}

\mon{May}

\parts{2}

\Date{MS received 1 January 2004}

\begin{abstract}
Let $R$ be a commutative noetherian ring and $f_{1}, \ldots, f_{r} \in R$.
In this article we give (cf. the Theorem in \S2) a criterion for
$f_{1}, \ldots, f_{r}$ to be regular sequence for a finitely generated
module over $R$ which strengthens and generalises a result in \cite{2}.
As an immediate consequence we deduce that if ${\rm V}(g_{1}, \ldots,
g_{r}) \subseteq {\rm V} (f_{1}, \ldots, f_{r})$ in Spec $R$ and if $f_{1}, \ldots,
f_{r}$ is a regular sequence in $R$, then $g_{1}, \ldots, g_{r}$ is
also a regular sequence in $R$.
\end{abstract}


\maketitle

\section{Regular sequences}

As there is no uniformity about the concept of regular sequence, we
first recall the following definitions that we shall use in this note.

\begin{definit}$\left.\right.$\vspace{.5pc}

\noindent {\rm Let $R$ be a commutative noetherian ring and $f_{1},
\ldots, f_{r} \in R$. We say that $f_{1}, \ldots, f_{r}$ is a {\it
strongly regular sequence} on a $R$-module $M$, if for every $i = 1,
\ldots, r$ the element $f_{i}$ is a non-zero divisor for $M/(f_{1},
\ldots, f_{i - 1})M$. The sequence $f_{1}, \ldots, f_{r}$ is called a
{\it regular sequence} on a $R$-module $M$, if for every ${\frak p} \in$
Supp($M/f_{1}, \ldots, f_{r})M)$, the sequence $f_{1}, \ldots, f_{r}$ in
the local ring $R_{\frak p}$ is a strongly regular sequence on the
$R_{\frak p}$-module $M_{\frak p}$.

Note that, in contrast to most of the standard text books, we do not
assume the $M \neq (f_{1}, \ldots, f_{r})M$ for a strongly regular
sequence $f_{1}, \ldots, f_{r}$. For general notations in commutative
algebra we also refer to \cite{1}.}\vspace{.5pc}
\end{definit}

If the sequence $f_{1}, \ldots, f_{r}$ is strongly regular respectively
regular on the $R$-module $M$, then the same is true for the sequence
$f_{1} \cdot 1_{S}, \ldots, f_{r}\cdot 1_{S}$ on the $S$-module
$S \otimes_{R} M$, where $S$ is an arbitrary flat noetherian
$R$-algebra.

Note that every sequence is a strongly regular as well as regular
sequence on the zero module. Further, it is clear that a strongly
regular sequence is a regular sequence but not conversely. For example:

\begin{exam}
Let $P : = k[X, Y, Z]$ be the polynomial ring in three indeterminates
over a field $k, {\frak p} : = P(X-1) + PZ, {\frak q} := PY$ and let $R
:= P/{\frak p} \cap {\frak q} = P/PY (X - 1) + PYZ$. Then $Z, X$ is a
regular sequence on the $P$-module $R$ but not a strongly regular
sequence.\vspace{.5pc}
\end{exam}

The difference between regular and strongly regular sequences is
well-illustrated in the following statement given in Chapter~II, 6.1 of
\cite{4}.

\begin{pro}$\left.\right.$\vspace{.5pc}

\noindent {\it Let $M$ be a finitely generated module over a noetherian
ring $R$ and let $f_{1}, \ldots, f_{r} \in R$. Then the following
conditions are equivalent{\rm :}
\begin{enumerate}
\renewcommand\labelenumi{\rm (\roman{enumi})}
\leftskip .1pc
\item $f_{1}, \ldots, f_{r}$ is a strongly regular sequence on $M$.\vspace{.15pc}

\item For every $s = 1, \ldots, r$ the sequence $f_{1}, \ldots, f_{s}$
is a regular sequence on $M$.
\end{enumerate}}\vspace{-.2pc}
\end{pro}

It can be easily seen that (see the proof of Proposition~3, Chapter~IV,
A, \S1 of \cite{5}) a sequence $f_{1}, \ldots, f_{r}$ in a commutative
noetherian ring $R$ is a regular sequence for a finitely generated
$R$-module $M$ if and only if the Koszul complex $K_{\bullet}\,(f_{1},
\ldots, f_{r}; M)$ gives a resolution of $M/(f_{1}, \ldots, f_{r})M$. In
particular, if $f_{1}, \ldots, f_{r}$ is a regular sequence on $M$, then
for every permutation $\sigma \in {\frak S}_{r}$ the sequence $f_{\sigma
1}, \ldots, f_{\sigma r}$ is also regular for $M$. Further, the above
proposition implies that the sequence $f_{\sigma 1}, \ldots, f_{\sigma
r}$ is strongly regular on $M$ for every $\sigma \in {\frak S}_{r}$ if
and only if all subsequences of $f_{1}, \ldots, f_{r}$ are regular on
$M$. For the sake of completeness let us recall Definition~2.

\begin{definit}$\left.\right.$\vspace{.5pc}

\noindent {\rm Let $(R, {\frak m}_{R})$ be a noetherian local ring and
let $M$ be a non-zero $R$-module. Then the length of a maximal regular
sequence on $M$ in the maximal ideal ${\frak m}_{R}$ is called the depth
of $M$ over $R$ and is denoted by depth$_{R} (M)$.}\vspace{.5pc}
\end{definit}

If $M$ is finitely generated then depth can be (cf. \cite{5},
Proposition and Definition~3, Chapter~IV, A, \S2) characterized by
\begin{equation*}
\hskip -4pc {(\ddagger)}\qquad\quad\ \ {\rm depth}_{R} (M) = \min \{ i \in
{\mathbb N} | {\rm Ext}_{R}^{i} (R/{\frak m}_{R}, M) \neq 0 \}.
\end{equation*}
A finitely generated $R$-module is called a {\it Cohen--Macaulay module}
if ${\rm dim}_{R}(M) = {\rm depth}_{R}(M)$.

\section{Theorem}

The following theorem is the main result of this note.

\begin{theo}
{\it Let $R$ be a commutative noetherian ring{\rm ,} $f_{1}, \ldots,
f_{r} \in R$ and let $M$ be a finitely generated $R$-module. Then the
following statements are equivalent{\rm :}

\begin{enumerate}
\renewcommand\labelenumi{\rm (\roman{enumi})}
\leftskip .15cm
\item $f_{1}, \ldots, f_{r}$ is a regular sequence on $M$.\vspace{.15pc}

\item  ${\rm depth}_{R_{\frak p}} (M_{\frak p}) \geq r$ for every
${\frak p} \in {\rm Supp} (M/(f_{1}, \ldots, f_{r})M)$.\vspace{.15pc}

\item  ${\rm depth}_{R_{\frak p}} (M_{\frak p}) \geq r$ for every
${\frak p} \in {\rm Ass} (M/(f_{1}, \ldots, f_{r})M)$.\vspace{-.5pc}
\end{enumerate}}
\end{theo}

\begin{proof}
The implications (i) $\Rightarrow$ (ii) $\Rightarrow$ (iii) are
trivial.\vspace{.5pc}

\noindent (ii) $\Rightarrow$ (i): We may assume that $R$ is local and
$f_{1}, \ldots, f_{r} \in {\frak m}_{R}$. Let ${\frak p} \in {\rm Ass}
(M)$ and let ${\frak q}$ be a minimal prime ideal in ${\rm V} ({\frak p}
+ Rf_{1} + \cdots + Rf_{r})$. Then ${\frak q} \in {\rm Supp}(M/(f_{1},
\ldots, f_{r})M) = {\rm Supp} (M) \cap {\rm V} (f_{1}, \ldots, f_{r})$
and so depth$_{R_{{\frak q}}} M_{\frak q} \geq r$ by (ii). Since ${\frak
p} \in {\rm Ass}(M)$, we have Hom$_{R_{\frak p}} (k ({\frak p}),
M_{\frak p}) \neq 0$ and so ${\rm Ext}_{R_{\frak q}}^{h} (k ({\frak q}),
M_{\frak q}) \neq 0$ by Chapter~6, \S18, Lemma~4 of \cite{3}, where $h :
= ht_{R/{\frak p}} ({\frak q}/{\frak p})$. Therefore $r \leq {\rm
depth}_{R_{\frak q}} M_{\frak q} \leq h$ (see $(\ddagger)$ in \S1). But
then $f_{1} \notin {\frak p}$, since otherwise $h \leq r -1$ by the
(generalised) Krull's theorem (see \cite{5}, Corollary~4, Chapter~III,
B, \S2). This proves that $f_{1}$ is a non-zero divisor for $M$. Now,
induction on $r$ completes the proof.

The implication (iii) $\Rightarrow$ (i) is proved in the lemma which is
given below. (In the proof of the lemma we use the implication (ii)
$\Rightarrow$ (i).)\hfill $\BBox$
\end{proof}

\begin{coro}{\rm (\cite{2}, Corollary~1)}$\left.\right.$\vspace{.5pc}

\noindent Let $R$ be a commutative noetherian ring{\rm ,} $f_{1},
\ldots, f_{r} \in R$ and let $M$ be a finitely generated $R$-module.
Then $f_{1}, \ldots, f_{r}$ is a regular sequence on $M$ if and only if
$f_{1}, \ldots, f_{r}$ is a regular sequence on $M_{\frak p}$ for every
${\frak p} \in {\rm Ass} (M/f_{1}, \ldots, f_{r})M)$.
\end{coro}

\begin{coro}$\left.\right.$\vspace{.5pc}

\noindent Let $R$ be a commutative noetherian ring and let \hbox{$f_{1},
\ldots, f_{r}, g_{1}, \ldots, g_{r} \in R$.} Let $M$ be a finitely
generated $R$-module such that ${\rm Supp} (M/(g_{1}, \ldots,
g_{r})M) \subseteq {\rm Supp}(M/(f_{1}, \ldots,$ $f_{r})M)$. Suppose that
$f_{1}, \ldots, f_{r}$ is a regular sequence on $M$. Then $g_{1},
\ldots, g_{r}$ is also a regular sequence on $M$. In particular{\rm ,} if
$\,{\rm V} (g_{1}, \ldots, g_{r}) \subseteq {\rm V} (f_{1}, \ldots,
f_{r})$ and if $f_{1}, \ldots, f_{r}$ is a regular sequence in $R${\rm ,}
then $g_{1}, \ldots, g_{r}$ is also a regular sequence in $R$.
\end{coro}

From the above equivalence we can also deduce the following well-known fact:

\begin{coro}{\rm (cf. \cite{5}, Theorem~2, Chapter~IV, B, \S2)}$\left.\right.$\vspace{.5pc}

\noindent If $M$ is a finitely generated Cohen--Macaulay module over a
noetherian local ring $R${\rm ,} then every system of parameters of $M$
is a regular sequence on $M$. In particular{\rm ,} in a Cohen--Macaulay
local ring every system of parameters is a regular sequence.
\end{coro}

Finally, we give a proof of the lemma which we have already used for the
proof of the implication (iii) $\Rightarrow$ (i) of the theorem.

\begin{lem}{\it
Let $R$ be a commutative noetherian ring{\rm ,} \hbox{$f_{1}, \ldots,
f_{r} \in R$} and let $M$ be a finitely generated $R$-module. Suppose
that \hbox{${\rm depth}_{R_{\frak p}} (M_{\frak p}) \geq r$} for every ${\frak
p} \in {\rm Ass} (M/(f_{1}, \ldots,$ $f_{r})M)$. Then $f_{1}, \ldots,
f_{r}$ is a regular sequence on $M$.}
\end{lem}

\begin{proof}
We shall prove by induction on $r$ the following implication:

\begin{enumerate}
\renewcommand\labelenumi{$(*)_{r}$:}
\leftskip .3cm
\item If ${\rm depth}_{R_{\frak p}} (M_{\frak p}) \geq r$ for every
${\frak p} \in {\rm Ass}(M/(f_{1}, \ldots, f_{r})M)$, then $f_{1},
\ldots, f_{r}$ is a regular sequence on $M$.
\end{enumerate}\vspace{-1pc}

\begin{prof}
Put $f := f_{1}$ and suppose that depth$_{R_{\frak p}} (M_{\frak p})
\geq 1$ for every ${\frak p} \in {\rm Ass} (M / fM)$. Then ${\rm Ass}
(M) \cap {\rm Ass} (M / fM) = \emptyset$. We shall show that $f$ is a non-zero
divisor for $M$. Suppose on the contrary that $f$ is a zero divisor on
$M$. By localising at a minimal prime ideal in ${\rm Ass} (M) \cap
{\rm V}(R f)$, we may assume that $R$ is a local ring,
depth$_{R} (M) = 0$ and that ${\rm Ass}(M) = \{ {\frak p}_{1}, \ldots,
{\frak p}_{m}, {\frak m}_{R} \}$ with ${\frak p}_{i} \notin {\rm V}
(Rf)$ for all $i = 1, \ldots, m$. Then $m \geq 1$. Let $Q_{1}, \ldots,
Q_{m}$ and $Q$ be the primary components corresponding to ${\frak
p}_{1}, \ldots, {\frak p}_{m}$ and ${\frak m}_{R}$ respectively and let
$0 = Q_{1} \cap \cdots \cap Q_{m} \cap Q$ be an irredundant primary
decomposition of the zero module in $M$. Let $N := Q_{1} \cap \cdots
\cap Q_{m}$. Then $N \neq 0, {\rm Ass}(M/N) = \{ {\frak p}_{1}, \ldots,
{\frak p}_{m} \}$ and $f$ is a non-zero divisor for $M/N$, since $f
\notin {\frak p}_{i}$ for all $i = 1, \ldots, m$. This implies that the
canonical homomorphism $N/fN \lr M/fM$ is injective. Further, since $Q$
is ${\frak m}_{R}$-primary in $M$, we have ${\frak m}_{R}^{n} N
\subseteq N \cap {\frak m}_{R}^{n} M \subseteq N \cap Q = 0$ for some $n
\in {\mathbb N}^{+}$, and hence $N$ has finite length. Therefore $N /
fN$ has finite length. But ${\rm depth}_{R}(M / fM) \geq 1$, since
${\frak m}_{R} \notin {\rm Ass}(M / fM)$ and therefore cannot contain any
submodules of finite length. This proves that $N / fN = 0$ and then $N =
0$ by Nakayama's lemma, which contradicts\break $N \neq 0$.
\end{prof}

\begin{prff}
We may assume that $R$ is local, $f_{1}, \ldots, f_{r + 1} \in {\frak
m}_{R}$ and $M \neq 0$. Now, we shall prove this implication by
induction on dim $(R)$. Clearly the induction starts at dim $(R) = 0$.
Put $\overline{M}_{r} := M / (f_{1}, \ldots, f_{r})M$ and
$\overline{M}_{r + 1} : = M / (f_{1}, \ldots, f_{r + 1})M$. Then by
induction hypothesis.
\begin{equation*}
\hskip -4pc (\dagger) \quad\qquad\ \ \hbox{$f_{1}, \ldots, f_{r + 1}\hbox{\ is
a regular sequence on\ }M_{\frak p}\hbox{\ for every\ } {\frak p} \in
{\rm Supp} (\overline{M}_{r + 1}) \backslash \{ {\frak m}_{R} \}.$}
\end{equation*}

$\left.\right.$\vspace{-1.5pc}

\noindent In particular, we have:
\begin{equation*}
\hskip -4pc (\dagger\dagger) \qquad\quad {\rm depth}_{R_{\frak p}}
(M_{\frak p}) \geq r + 1 \ \hbox{ for every }\ {\frak p} \in {\rm Supp}
(\overline{M}_{r + 1})\backslash \{ {\frak m}_{R} \}.
\end{equation*}
\end{prff}
We consider two cases:

\begin{case}\hskip -.3pc {\rm \hbox{${\frak m}_{R} \in {\rm Ass}
(\overline{M}_{r + 1}).$}}\ \ {\rm In this case, by assumption in \hbox{$(*)_{r
+ 1}, {\rm depth}_{R} (M) \geq r + 1$.} Now, use (ii) $\Rightarrow$ (i)
of the theorem to conclude that $f_{1}, \ldots, f_{r + 1}$ is a regular
sequence on~$M$.}
\end{case}

\begin{case}\hskip -.3pc {\rm \hbox{${\frak m}_{R} \notin {\rm
Ass}(\overline{M}_{r + 1}).$}}\ \ {\rm In this case ${\rm Ass}({\overline
M}_{r}) \ \cap \ {\rm Ass}({\overline M}_{r + 1}) = \emptyset$, since
${\rm depth}_{R_{\frak p}}(\overline{M}_{r})_{\frak p} \geq 1$ for every
${\frak p} \in {\rm Ass}(\overline{M}_{r + 1}) \backslash \{ {\frak
m}_{R} \}$ by $(\dagger \dagger)$. Therefore by $(*)_{1}, f_{r + 1}$ is
a non-zero divisor on $\overline{M}_{r}$. Now, it remains to show that
the sequence $f_{1}, \ldots, f_{r}$ is a regular sequence on $M$. For
this, let ${\frak p} \in {\rm Ass} (\overline{M}_{r})$. Since $f_{r +
1}$ is a non-zero divisor for $\overline{M}_{r}$, there exists ${\frak
q} \in {\rm Ass}(\overline{M}_{r + 1})$ such that ${\frak p} \subseteq
{\frak q}$. Note that ${\frak q} \neq {\frak m}_{R}$ and that $f_{1},
\ldots, f_{r}$ is a regular sequence on $M_{\frak q}$ by $(\dagger)$ and
hence in particular for $M_{\frak p}$. This proves that depth$_{R_{\frak
p}} (M_{\frak p}) \geq r$ for every ${\frak p} \in {\rm
Ass}(\overline{M}_{r})$ and hence $f_{1}, \ldots, f_{r}$ is a regular
sequence on $M$ by $(*)_{r}$.}\hfill $\BBox\ $
\end{case}\vspace{-1pc}
\end{proof}

\section*{Acknowledgements}

Part of this work was done while the first author was visiting Germany
during April--June~2001 under a grant from DAAD, Germany.\ \,The first
author thanks DAAD, Germany for financial support. The authors sincerely
thank Harmut Wiebe for stimulating discussions.

\end{document}